\documentclass[a4paper,11pt,twoside,reqno]{amsart}

\usepackage[utf8]{inputenc}
\usepackage[plainpages=false,pdfpagelabels=true]{hyperref}
\usepackage{amssymb,amsthm}
\usepackage[margin=1in]{geometry}
\usepackage{slashed}

\newtheorem{Satz}{Theorem}[section]

\newtheorem{Lem}[Satz]{Lemma}

\newtheorem{Cor}[Satz]{Corollary}
\theoremstyle{definition}

\newtheorem{Bem}[Satz]{Remark}

\newcommand{\tr}{\operatorname{Tr}}

\allowdisplaybreaks[1]

\renewcommand{\epsilon}{\varepsilon}

\newcommand{\R}{\ensuremath{\mathbb{R}}}

\numberwithin{equation}{section}

\title{On harmonic and biharmonic maps from gradient Ricci solitons}
\author{Volker Branding}
\date{\today}
\address{University of Vienna, Faculty of Mathematics\\
Oskar-Morgenstern-Platz 1, 1090 Vienna, Austria\\}
\email{volker.branding@univie.ac.at}
\thanks{The author gratefully acknowledges the support of the Austrian Science Fund (FWF) through the START-project (Y963) of Michael Eichmair and 
the project "Geometric Analysis of Biwave Maps" (P 34853).
}
\begin{document}

\begin{abstract}
We study harmonic and biharmonic maps from gradient Ricci solitons.
We derive a number of analytic and geometric conditions under which harmonic
maps are constant and which force biharmonic maps to be harmonic.
In particular, we show that biharmonic maps of finite energy from
the two-dimensional cigar soliton must be harmonic.
\end{abstract} 

\maketitle

\section{Introduction and results}
One of the central aims in the geometric calculations of variations is to
find interesting maps between Riemannian manifolds. 
This can often be achieved by extremizing a given energy functional.
For a concrete energy functional under consideration one aims
to obtain a classification of its critical points, i.e. one wants to
know under which conditions critical points can exist or are obstructed.
In this article, we will focus on harmonic and biharmonic maps
between Riemannian manifolds and
establish a number of classification results for these.

Becoming more technical, we assume that \((M,g)\) and \((N,h)\) are two Riemannian manifolds.
Consider a map \(\phi\colon M\to N\), its energy is defined by
\begin{equation}
\label{eq:energy}
E_1(\phi)=E(\phi)=\int_M|d\phi|^2dv_g.
\end{equation}
It is well-known that the critical points of \eqref{eq:energy} are \emph{harmonic maps}.
They are characterized by the vanishing of the \emph{tension field}
which is defined as follows
\begin{equation}
\tau(\phi)=\tr_g\bar\nabla d\phi,
\end{equation}
where \(\bar\nabla\) denotes the connection on the vector bundle \(\phi^\ast TN\).
The harmonic map equation \(\tau(\phi)=0\) is a semilinear elliptic partial
differential equation of second order, for more details we refer to the book \cite{MR1391729}.
Of course, the harmonic map equation admits trivial solutions,
namely in the case that the map \(\phi\) maps to a point \(q\in N\).

A higher order generalization of harmonic maps that receives growing attention
are the so-called \emph{biharmonic maps}.
These are critical points of the bienergy functional
\begin{equation}
\label{eq:bienergy}
E_2(\phi)=\int_M|\tau(\phi)|^2dv_g
\end{equation}
and are characterized by the vanishing of the so-called
\emph{bitension field} which is explicitly given by
\begin{equation}
\label{eq:bitension}
\tau_2(\phi):=\bar\Delta\tau(\phi)+\tr_gR^N(\tau(\phi),d\phi)d\phi.
\end{equation}
Here, \(\bar\Delta\) denotes the connection Laplacian on the vector bundle 
\(\phi^\ast TN\) for which we choose the analysts sign convention, that is
\begin{equation*}
\bar\Delta:=\tr_g(\bar\nabla\bar\nabla-\bar\nabla_{\nabla}).
\end{equation*}
In contrast to the harmonic map equation, the biharmonic map equation 
\(\tau_2(\phi)=0\),
is a fourth order semilinear elliptic partial differential equation which leads to additional technical difficulties. In particular, tools like the maximum principle,
which are well-adapted to partial differential equations of second order,
are no longer applicable in full generality in the case of biharmonic maps.

For an overview on biharmonic maps we refer to the survey article \cite{MR2301373}
and the recent book \cite{MR4265170}.
Both harmonic and biharmonic maps have important applications in analysis and geometry,
but are also utilized in elasticity theory and quantum field theory.

For an overview on further higher order variational problems in Riemannian geometry, generalizing harmonic and biharmonic maps,
we refer to the paper \cite{MR4106647}.

It can be directly seen that every harmonic map, that is a solution of \(\tau(\phi)=0\),
automatically solves the equation for biharmonic maps which is \(\tau_2(\phi)=0\).
Hence, one is interested in constructing biharmonic maps which are non-harmonic,
the latter are called \emph{proper biharmonic}.
On the other hand, under certain geometric and analytic assumptions a biharmonic map 
necessarily needs to be harmonic.
For example, if both \(M,N\) are compact and \(N\) has negative sectional curvature
then a direct application of the maximum principle shows that every biharmonic map is harmonic.
More results of this kind can be found in  \cite{MR3834926,MR4040175}.

In this article we will derive a number of inequalities that allow to get an understanding when harmonic maps are constant and
when biharmonic maps need to be harmonic assuming that the domain manifold \(M\)
is equipped with a special kind of Riemannian metric which are gradient Ricci solitons.

In order to approach these inequalities let us first recall a number of results on 
gradient Ricci solitons.
The concept of Ricci solitons is of great importance in the study of the Ricci flow
as Ricci solitons evolve just by diffeomorphisms and homotheties of the initial metric 
under Ricci flow. For an introduction to Ricci flow we refer to the book \cite{MR2265040}.
A special class of Ricci solitons are the so-called \emph{gradient Ricci solitons}.
They are characterized by the following equation:
\begin{equation}
\label{eq:soliton}
\operatorname{Ric^M}+\nabla^2f=\lambda g
\end{equation}
Here, \(\operatorname{Ric^M}\) represents the Ricci curvature of the Riemannian 
manifold \((M,g)\), \(\nabla^2f\) is the Hessian of a function \(f\in C^\infty(M)\) and \(\lambda\in\R\).

A gradient Ricci soliton is called \emph{steady} if \(\lambda =0\),
\emph{shrinking} if \(\lambda>0\) and \emph{expanding} if \(\lambda <0\).

In the case of a two-dimensional manifold gradient Ricci solitons can be completely
classified, see \cite{MR3340295,ram}. In higher dimensions it is substantially more difficult to obtain such a classification, see for example \cite{MR2507581}, \cite{MR2818729} for important results on this matter.

A gradient Ricci soliton is called \emph{trivial} if \(f=const\).
In the case of \(M\) being compact gradient Ricci solitons are most often
trivial as can be seen by taking the trace of \eqref{eq:soliton} and applying
the maximum principle. For this reason, we only consider the case of a non-compact manifold \(M\) in this manuscript.

Let us mention several results from the literature which are closely
connected to the content of this article.
Sealey showed that harmonic maps of finite energy from Euclidean and hyperbolic
space of dimensions bigger than two must be constant \cite{MR654088}.
Further results of this kind based on monotonicity formulas
are presented in \cite[Section 2.3]{MR1391729}.
Harmonic functions on gradient Ricci solitons were studied in \cite{MR3023848}
by Munteanu and Sesum.
Rimoldi investigated \(f\)-harmonic maps, which are critical points of a weighted
version of harmonic maps, and their relation to gradient Ricci solitons in
\cite{MR2818729}.
Moreover, \(L^2\)-harmonic forms on gradient shrinking Ricci solitons were investigated by Yun in \cite{MR3668864}.

A class of solitons closely related to this article are the so-called
\emph{Ricci-harmonic solitons} studied in \cite{MR3406677}.
These solitons arise in the context of the Ricci-harmonic flow
which is a combination of both Ricci and harmonic map heat flow
introduced in \cite{MR2961788}.

Let us briefly describe the strategy that we are using in order to establish our main results. For both harmonic and biharmonic maps there exists an associated stress-energy tensor which is divergence free. Testing this conservation law with the equation
for a gradient Ricci soliton, using a cutoff-function in order to be able 
to employ integration by parts, we are led to a number of energy inequalities
from which we can deduce various vanishing results.

In the following we set \(m:=\dim M\) and present the main results of this article.
\begin{Satz}
\label{thm:liouville-harmonic}
Assume that \((M,g,f)\) is a complete non-compact gradient Ricci soliton. 
Let \(\phi\colon M\to N\) be a smooth harmonic map which satisfies
\begin{equation}
\label{eq:harmonic-finite-energy}
\lim_{R\to\infty}\frac{1}{R}\int_{B_R}|\nabla f||d\phi|^2d\mu=0,
\end{equation}
where \(B_R\) denotes the geodesic ball of radius \(R\).
Then the following inequality holds
\begin{equation}
\label{inequality-harmonic}
\int_M\big(\lambda(m-2)-\operatorname{Scal}^M\big)|d\phi|^2dv_g
+2\int_M\langle d\phi(\operatorname{Ric}^M(e_i)),d\phi(e_i)\rangle dv_g\leq 0,
\end{equation}
where \(\operatorname{Scal}^M\) represents the scalar curvature of the manifold \(M\)
and \(\{e_i\},i=1,\ldots,m\) denotes an orthonormal frame field on \(M\).
\end{Satz}
In order to obtain any information from \eqref{inequality-harmonic}
one of course has to make additional curvature assumptions
that guarantee that the left-hand side of \eqref{inequality-harmonic} is positive.
One possibility to obtain a kind of Liouville theorem is given by the following
\begin{Cor} 
\label{cor:harmonic}
Assume that \((M,g,f)\) is a complete, non-compact gradient Ricci soliton
with \(m>2\) and 
let \(\phi\colon M\to N\) be a smooth harmonic map which satisfies \eqref{eq:harmonic-finite-energy}.
If 
\begin{equation*}
\int_M\langle d\phi(\operatorname{Ric}^M(e_i)),d\phi(e_i)\rangle dv_g\geq 0,
\end{equation*}
and \(\lambda(m-2)-\operatorname{Scal}^M>0\)
then \(\phi\) is a constant map.
\end{Cor}

\begin{Bem}
\begin{enumerate}
 \item 
Note that the inequality \eqref{inequality-harmonic} does not contain any information if
\(m=2\) as we have that \(2\operatorname{Ric}^M=\operatorname{Scal}^Mg\).
\item The assumptions in Corollary \ref{cor:harmonic} seem to be restrictive
but there does not seem to be any kind of contradiction.
\item Note that Corollary \ref{cor:harmonic}
recovers the result of Sealey \cite{MR654088} 
in the case of \((M,g)=(\R^m,\delta)\),
where \(\delta\) represents the flat Euclidean metric, by choosing 
\(f=\frac{|x|^2}{2}\). This particular choice of \(f\) is 
called \emph{Gaussian soliton} in the literature.
\item
It is well-known that a harmonic map of finite energy from a complete, non-compact manifold of positive Ricci curvature to a Riemannian manifold of negative sectional curvature must be constant due to a celebrated result of Schoen and Yau \cite{MR438388}.
The previous results are in the same spirit but here we only
make assumptions on the geometry of the domain and do not require the target
to have negative curvature.
\end{enumerate}
\end{Bem}

In addition, making use of the Bochner formula for harmonic maps, we will also
establish the following statement extending a result for 
harmonic functions \cite[Theorem 4.1]{MR3023848}:

\begin{Satz}
\label{thm:liouville-harmonic-curvature}
Let \((M,g,f)\) be a complete, non-compact steady gradient Ricci soliton of
infinite volume and \((N,h)\)
a Riemannian manifold of non-positive sectional curvature.
In addition, let \(\phi\colon M\to N\) be a harmonic map satisfying
\begin{equation}
\int_M|d\phi|^2dv_g<\infty.
\end{equation}
Then \(\phi\) must be a constant map.
\end{Satz}

Besides the above Liouville-type results for harmonic maps from gradient Ricci solitons
we will also prove the following corresponding results for biharmonic maps.

\begin{Satz}
\label{thm:liouville-biharmonic}
Assume that \((M,g,f)\) is a complete, non-compact gradient Ricci soliton
with \(\operatorname{Ric}^M\leq C\) and \(|\nabla f|<\infty\).
Let \(\phi\colon M\to N\) be a smooth biharmonic map with finite energy, that is
\begin{equation}
\label{eq:biharmonic-finite-energy}
\int_M(|d\phi|^2+|\bar\nabla d\phi|^2)dv_g<\infty.
\end{equation}
Then the following inequality holds
\begin{equation}
\label{identity-biharmonic}
\int_M\big(\lambda(m-4)-\operatorname{Scal}^M\big)|\tau(\phi)|^2dv_g
+4\int_M\operatorname{Ric}^M(e_i,e_j)\langle\bar\nabla d\phi(e_i,e_j),\tau(\phi)\rangle dv_g
\leq 0.
\end{equation}
\end{Satz}

\begin{Bem}
\begin{enumerate}
 \item It is clear that the inequality \eqref{identity-biharmonic} has a similar structure
as the corresponding inequality for harmonic maps \eqref{inequality-harmonic}.
While \eqref{identity-biharmonic} reflects that biharmonic maps are a fourth order equation, \eqref{inequality-harmonic} clearly shows the second order character of harmonic maps. However note that in Theorem \ref{thm:liouville-biharmonic} we need to require an upper bound on the Ricci curvature which is not necessary in the case
of harmonic maps (Theorem \ref{thm:liouville-harmonic}).
Moreover, note that the finite energy assumption \eqref{eq:biharmonic-finite-energy}
involves the full second covariant derivative \(\bar\nabla d\phi\) and not only
its trace given by the tension field \(\tau(\phi)\).
\item If one tries to derive corresponding energy inequalities for polyharmonic maps,
then one can expect that it is necessary to also require upper bounds on the covariant
derivatives of the Ricci curvature.
\end{enumerate}
\end{Bem}

In addition, we can also give
the following kind of Liouville-type result obtained from Theorem \ref{thm:liouville-biharmonic}.

\begin{Cor}
\label{cor:biharmonic}
Assume that \((M,g,f)\) is a complete, non-compact gradient Ricci soliton
with \(\operatorname{Ric}^M\leq C\) and \(|\nabla f|<\infty\).
Let \(\phi\colon M\to N\) be a smooth biharmonic map with finite energy.
If
\begin{equation*}
\int_M\operatorname{Ric}^M(e_i,e_j)\langle\bar\nabla d\phi(e_i,e_j),\tau(\phi)\rangle dv_g
\geq 0
\end{equation*}
and \(\lambda(m-4)-\operatorname{Scal}^M>0\)
then \(\phi\) must be a harmonic map.
\end{Cor}

Let us recall a famous example of a steady gradient Ricci soliton,
which is the so-called \emph{Hamilton's cigar soliton} in two dimensions.
In the physics literature it is known as \emph{Witten's black hole},
see \cite[p. 10]{MR2265040} for some more details.
It is given by the following data
\begin{equation}
\label{dfn:cigar-soliton}
(M,g,f)=\big(\R^2,\frac{dx^2+dy^2}{1+x^2+y^2},-\log(1+x^2+y^2)\big)
\end{equation}
and has positive scalar curvature \(\operatorname{Scal}^M=\frac{1}{1+x^2+y^2}>0\).

For biharmonic maps from the cigar soliton we can give the following

\begin{Cor}
\label{cor:cigar-biharmonic}
Assume that \((M,g,f)\) is the two-dimensional cigar soliton defined by 
\eqref{dfn:cigar-soliton}.
Let \(\phi\colon M\to N\) be a smooth biharmonic map with finite energy.
Then \(\phi\) is harmonic.
\end{Cor}

\begin{Bem}
\begin{enumerate}
\item In the case that \((M,g,f)\) is a steady gradient Ricci soliton,
that is \(\lambda=0\) we do not need to require that \(|\nabla f|<\infty\)
as this condition is automatically satisfied, see Section 2.

However, in general, we have to make the assumption that \(|\nabla f|<\infty\)
although very often this condition may not be necessary but 
this of course depends on the concrete gradient Ricci soliton.

If we want to drop the assumption \(|\nabla f|<\infty\) in Theorem \ref{thm:liouville-biharmonic}, then in addition to
the finite energy assumption \eqref{eq:biharmonic-finite-energy}, we have to demand
\begin{equation}
\lim_{R\to\infty}\frac{1}{R}\int_{B_R}|\nabla f|\big(|d\phi|^2+|\bar\nabla d\phi|^2\big)d\mu=0,
\end{equation}
where \(B_R\) denotes the geodesic ball of radius \(R\).

\item Of course there may be geometric configurations under which Corollary \ref{cor:biharmonic}
  states that a biharmonic map with finite energy has to be constant rather
  than just being harmonic.
\item Note that Theorem \ref{thm:liouville-biharmonic}
recovers a result of Baird et al. \cite[Theorem 3.4]{MR2604617} 
in the case of \((M,g)=(\R^m,\delta)\),
where \(\delta\) represents the flat Euclidean metric, by choosing 
\(f=\frac{|x|^2}{2}\), see also \cite[Theorem 3.2]{MR4142862} for a 
different method of proof. 

\item It does not seem to be possible to prove a result of the form of Theorem \ref{thm:liouville-harmonic-curvature} for biharmonic maps as the Ricci curvature of the domain does not enter in the Bochner
formula for biharmonic maps.
\item If we rescale the metric \(g\to r^2g\) for \(r\in\R\) in the equation for a gradient Ricci soliton \eqref{eq:soliton} then it changes to 
\begin{equation*}
\operatorname{Ric^M}+\nabla^2f=r^2\lambda g.
\end{equation*}
Hence, many of the geometric conditions that we can impose to 
deduce a vanishing result from Theorems \ref{thm:liouville-harmonic}, 
\ref{thm:liouville-biharmonic} may not be invariant under rescaling
of the metric.
\end{enumerate}
\end{Bem}

This article is organized as follows:
In Section 2 we recall the relevant background material on gradient Ricci solitons
and stress-energy tensors that is utilized in this article.
Afterwards, in Section 3, we provide the proofs of the main results.
Finally, Section 4 provides some remarks on harmonic and biharmonic maps
in the case that the domain manifold is a gradient Yamabe soliton.

Throughout this manuscript we make use of the summation convention, that is
we sum over repeated indices. We use the symbol \(\langle\cdot,\cdot\rangle\)
to denote various different scalar products.
Moreover, the letter \(C\) denotes a generic positive constant whose
value may change from line to line.

We use the following sign convention for the Riemannian curvature tensor
\[R(X,Y)Z=\nabla_X\nabla_YZ-\nabla_Y\nabla_XZ-\nabla_{[X,Y]}Z\]
for vector fields \(X,Y,Z\).
Concerning the Laplacian we use the analysts sign convention, such that
\(\Delta \xi=\xi''\) for \(\xi\in C^\infty(\R)\).

\par\medskip
\textbf{Acknowledgements:}
The author would like to thank the reviewers for many helpful comments 
which helped to substantially improve the content of the article.

\section{Some background material}
In this section we recall several well-known results on gradient Ricci solitons
and stress-energy tensors which are applied in the proofs of the main results.

We will often make use of the following identity
\begin{equation}
\label{eq:identity-curvature}
\operatorname{div}\operatorname{Ric}=\frac{1}{2}\nabla\operatorname{Scal}
\end{equation}
which follows from contracting the second Bianchi identity twice.
Note that \eqref{eq:identity-curvature} holds on every Riemannian manifold.

\subsection{Some properties of gradient Ricci solitons}
First, we recall the following facts on gradient Ricci soliton,
for more details we refer to \cite{MR2448435}.

\begin{Lem}
Assume that \((M,g,f)\) is a gradient Ricci soliton.
Then the following equation holds
\begin{equation}
\label{identity-soliton-const}
\operatorname{Scal}^M+|\nabla f|^2-2\lambda f=C
\end{equation}
for some constant \(C\).
\end{Lem}
\begin{proof}
First of all, we note that by taking the trace of \eqref{eq:soliton}
we get 
\begin{equation}
\Delta f=m\lambda-\operatorname{Scal}^M.
\end{equation}
By \(R_{ij}\) we denote the components of the Ricci tensor 
\(\operatorname{Ric}^M\).
Now, taking the divergence of \eqref{eq:soliton} we get
\begin{align*}
\operatorname{div}\operatorname{Ric}_i=&
\nabla^jR_{ij}\\
=&-\nabla^j\nabla_i\nabla_jf \\
=&-\nabla_i\Delta f-R_{ij}\nabla^jf \\
=&\nabla_i\operatorname{Scal}^M-R_{ij}\nabla^jf.
\end{align*}
Using \eqref{eq:identity-curvature} we deduce that
\begin{equation}
\label{eq:identity-a}
\nabla_i\operatorname{Scal}^M=2R_{ij}\nabla^jf.
\end{equation}
Combining \eqref{eq:identity-a} and \eqref{eq:soliton} we find
\begin{equation*}
\nabla_i\big(\operatorname{Scal}^M+|\nabla f|^2-2\lambda f\big)=0
\end{equation*}
completing the proof.
\end{proof}

\begin{Lem}
\label{lem:bound-nabla-f}
Assume that \((M,g,f)\) is a steady gradient Ricci soliton, then
\(\operatorname{Scal}^M\geq 0\) and 
\begin{equation}
\label{eq:steady-bound-f}
|\nabla f|^2\leq C
\end{equation}
for some positive constant \(C\).
\end{Lem}
\begin{proof}
The proof uses ideas from Ricci flow.
It is well-known that every steady Ricci soliton 
is an ancient solution to the Ricci flow.
However, by \cite[Corollary 2.5]{MR2520796} we know that any ancient
smooth complete solution to the Ricci flow must have non-negative
scalar curvature \(\operatorname{Scal}^M\geq 0\).
The bound on \(|\nabla f|^2\) now follows from \eqref{identity-soliton-const}.
\end{proof}

\subsection{Stress-energy tensors}
Besides the aforementioned results on gradient Ricci solitons we also recall
the \emph{stress-energy tensors} for both harmonic and biharmonic maps.
The stress-energy tensors arise by varying the energies \eqref{eq:energy} and \eqref{eq:bienergy} with respect to the metric on the domain.

In the case of harmonic maps the stress-energy tensor was first derived in \cite{MR655417}
and it is given by
\begin{equation}
\label{stress-energy-harmonic}
S_1(X,Y)=\frac{1}{2}|d\phi|^2g(X,Y)-\langle d\phi(X),d\phi(Y)\rangle,
\end{equation}
where \(X,Y\) are vector fields on \(M\).

A direct computation shows that the stress-energy tensor \eqref{stress-energy-harmonic}
satisfies the equation
\begin{equation}
\label{eq:conservation-harmonic}
\operatorname{div} S_1=-\langle\tau(\phi),d\phi\rangle.
\end{equation}
In particular, the stress-energy tensor \(S_1\) is conserved (has vanishing divergence)
if \(\phi\colon M\to N\) is a harmonic map, that is a solution of \(\tau(\phi)=0\).

Now, let us reconsider 
the stress-energy associated with the bienergy \eqref{eq:bienergy} which is given by
\begin{equation}
\label{stress-energy-bienery}
S_2(X,Y)=\big(\frac{1}{2}|\tau(\phi)|^2+\langle d\phi,\bar\nabla\tau(\phi)\rangle\big)g(X,Y)
-\langle d\phi(X),\bar\nabla_Y\tau(\phi)\rangle-\langle d\phi(Y),\bar\nabla_X\tau(\phi)\rangle,
\end{equation}
where \(X,Y\) are again vector fields on \(M\).
It was first stated by Jiang in \cite{MR891928} and later systematically studied 
by Loubeau, Montaldo and Oniciuc in \cite{MR2395125}.

The stress-energy tensor \eqref{stress-energy-bienery} satisfies the following conservation law
\begin{equation}
\label{eq:conservation-biharmonic}
\operatorname{div} S_2=-\langle\tau_2(\phi),d\phi\rangle.
\end{equation}
As in the case of harmonic maps, the stress-energy tensor \(S_2\) is conserved (has vanishing divergence)
if \(\phi\colon M\to N\) is a biharmonic map, that is a solution of \(\tau_2(\phi)=0\).
Recall that the bitension field \(\tau_2(\phi)\) is defined in \eqref{eq:bitension}.

The fact that the stress-energy tensors are conserved is a direct consequence of Noether's theorem as the energies \eqref{eq:energy} and \eqref{eq:bienergy} are
invariant under diffeomorphisms on the domain.

For more details on stress-energy tensors in the context of harmonic maps, 
we refer to \cite{MR4007262}
where the stress-energy tensor for polyharmonic maps was investigated in detail.

\section{Proof of the main results}
In this section we provide the proofs of the main results.
First, we will establish the following 
\begin{Lem}
Let \((M,g,f)\) be a complete, non-compact gradient Ricci soliton.
Suppose that \(\phi\colon M\to N\) is a smooth harmonic map and \(\eta\in C^\infty(M)\) with
compact support. Then the following identity holds
\begin{equation}
\label{identity-harmonic-a}
0=\frac{1}{2}\int_M\eta^2|d\phi|^2\Delta fdv_g
-\int_M\eta^2\nabla_{e_i}\nabla_{e_j}f\langle d\phi(e_i),d\phi(e_j)\rangle dv_g
+\int_M S_1(\nabla\eta^2,\nabla f)dv_g,
\end{equation}
where \(S_1\) represents the stress-energy tensor associated with harmonic maps.
\end{Lem}

\begin{proof}
By assumption the map \(\phi\) is harmonic hence the stress-energy tensor \(S_1\)
given by \eqref{stress-energy-harmonic} is divergence-free.
Let \(\{e_i\},i=1,\ldots,m\) be an orthonormal basis of \(TM\) that satisfies 
\(\nabla_{e_i}e_j=0,i,j=1,\ldots,m\) at a fixed point \(p\in M\).
As \(\eta\in C^\infty(M)\) has compact support we can use integration
by parts to calculate
\begin{align*}
0=&-\int_M\eta^2\langle\nabla_{e_i}f,\nabla_{e_j} S_1(e_i,e_j)\rangle dv_g \\
=&\int_M\eta^2\langle\nabla^2f,S_1\rangle dv_g
+\int_M S_1(\nabla\eta^2,\nabla f)dv_g 
\end{align*}
and making use of the definition of the stress-energy tensor \(S_1\) completes the proof.
\end{proof}

\begin{proof}[Proof of Theorem \ref{thm:liouville-harmonic}]
First of all, by taking the trace of \eqref{eq:soliton} we get
\begin{equation*}
\Delta f=m\lambda-\operatorname{Scal}^M. 
\end{equation*}
Now, we choose the function \(\eta\) as follows:
Let  \(0\leq\eta\leq 1\) on \(M\) be such that
\begin{equation*}
\eta(x)=1\textrm{ for } x\in B_R(x_0),\qquad \eta(x)=0\textrm{ for } x\in B_{2R}(x_0),\qquad |\nabla\eta|\leq\frac{C}{R}\textrm{ for } x\in M,
\end{equation*}
where \(B_R(x_0)\) denotes the geodesic ball around the point \(x_0\) with radius \(R\).

Using this identity, the defining equation of a gradient Ricci soliton \eqref{eq:soliton}
and \eqref{identity-harmonic-a} 
we obtain
\begin{equation*}
\int_M\eta^2\big(\lambda(m-2)-\operatorname{Scal}^M\big)|d\phi|^2dv_g
+2\int_M\eta^2\langle d\phi(\operatorname{Ric}^M(e_i)),d\phi(e_i)\rangle dv_g
=-2\int_M S_1(\nabla\eta^2,\nabla f)dv_g.
\end{equation*}
Now, we note that
\begin{align*}
\big|\int_M S_1(\nabla\eta^2,\nabla f)dv_g\big|&\leq
C\int_M|\eta||\nabla\eta||\nabla f||d\phi|^2dv_g\\
&\leq 
\frac{C}{R}\int_{B_R}|\nabla f||d\phi|^2d\mu.
\end{align*}
Due to the finite energy assumption \eqref{eq:harmonic-finite-energy}
this term will vanish as we take the limit \(R\to\infty\).
The proof is now complete.
\end{proof}

\begin{proof}[Proof of Theorem \ref{thm:liouville-harmonic-curvature}]  
First, recall that in the case of a steady gradient Ricci soliton,
that is \(\lambda=0\), we have
\begin{equation}
\operatorname{Ric}^M=-\operatorname{Hess}f,\qquad \operatorname{Scal}^M=-\Delta f.
\end{equation}
Inserting these identities into \eqref{identity-harmonic-a} we find
\begin{equation*}
0=\frac{1}{2}\int_M\eta^2|d\phi|^2\operatorname{Scal}^Mdv_g
-\int_M\eta^2\langle d\phi(\operatorname{Ric}^M(e_i)),d\phi(e_i)\rangle dv_g
-\int_M S_1(\nabla\eta^2,\nabla f)dv_g,
\end{equation*}
where \(S_1\) is the stress-energy tensor associated with harmonic maps \eqref{stress-energy-harmonic}.

Now, making use of the Bochner formula for harmonic maps,
see for example \cite[Prop. 1.3.5]{MR1391729}, we find
\begin{align*}
\Delta\frac{1}{2}|d\phi|^2=&|\bar\nabla d\phi|^2+\langle d\phi(\operatorname{Ric}^M(e_i)),d\phi(e_i)\rangle
-\langle R^N(d\phi(e_i),d\phi(e_j))d\phi(e_j),d\phi(e_i)\rangle \\
\nonumber\geq&\big|\nabla |d\phi|\big|^2+\langle d\phi(\operatorname{Ric}^M(e_i)),d\phi(e_i)\rangle,
\end{align*}
where we used the assumption of non-positive sectional curvature of the manifold \(N\)
and the Kato inequality
in the second step.
Combining the previous equations yields
\begin{align*}
\frac{1}{2}\int_M\eta^2|d\phi|^2\operatorname{Scal}^Mdv_g+\int_M\eta^2\big|\nabla |d\phi|\big|^2dv_g&\leq&
\frac{1}{2}\int_M\eta^2\Delta|d\phi|^2dv_g+\int_M S_1(\nabla\eta^2,\nabla f)dv_g.
\end{align*}

As in the previous proof we choose \(\eta\) as follows:
Let  \(0\leq\eta\leq 1\) on \(M\) be such that
\begin{equation*}
\eta(x)=1\textrm{ for } x\in B_R(x_0),\qquad \eta(x)=0\textrm{ for } x\in B_{2R}(x_0),\qquad |\nabla\eta|\leq\frac{C}{R}\textrm{ for } x\in M,
\end{equation*}
where \(B_R(x_0)\) denotes the geodesic ball around the point \(x_0\) with radius \(R\).
Now, we note that 
\begin{equation*}
\frac{1}{2}\int_M\eta^2\Delta|d\phi|^2dv_g=-\int_M\eta\nabla\eta\nabla|d\phi|^2dv_g
=-2\int_M\eta\nabla\eta|d\phi|\nabla|d\phi|dv_g.
\end{equation*}
Using Young's inequality this yields
\begin{equation*}
\frac{1}{2}\int_M\eta^2\Delta|d\phi|^2dv_g\leq
\frac{1}{2}\int_M\eta^2\big|\nabla |d\phi|\big|^2dv_g
+\frac{C}{R^2}\int_M|d\phi|^2dv_g,
\end{equation*}
where we also employed the properties of the cut-off function \(\eta\).

Recall that we
have
\begin{equation*}
\big|\int_M S_1(\nabla\eta^2,\nabla f)dv_g\big|\leq
\frac{C}{R}\int_M|\nabla f||d\phi|^2dv_g.
\end{equation*}
By assumption \((M,g,f)\) is a steady gradient Ricci soliton
such that \(|\nabla f|\leq C\) due to \eqref{eq:steady-bound-f}.

Hence, by combining the previous identities, we find
\begin{equation*}
\frac{1}{2}\int_M\eta^2|d\phi|^2\operatorname{Scal}^Mdv_g+\frac{1}{2}\int_M\eta^2\big|\nabla |d\phi|\big|^2dv_g\leq
\frac{C}{R^2}\int_M|d\phi|^2dv_g.
\end{equation*}
Now, letting \(R\to\infty\) and using the assumption of finite energy \eqref{eq:harmonic-finite-energy} we get
\begin{equation*}
\int_M\eta^2|d\phi|^2\operatorname{Scal}^Mdv_g+\int_M\eta^2\big|\nabla |d\phi|\big|^2dv_g\leq 0.
\end{equation*}
Hence, we may conclude that \(|d\phi|^2=C\) for a positive constant \(C\).
However, we require that the manifold \((M,g,f)\) has infinite volume
and due to the finite energy assumption we conclude that \(|d\phi|^2=0\).
\end{proof}

Before we turn to the proof of Theorem \ref{thm:liouville-biharmonic} we establish
a technical lemma.

\begin{Lem}
Assume that \((M,g,f)\) is a gradient Ricci soliton and
let \(\phi\colon M\to N\) be a smooth biharmonic map.
For \(\eta\in C^\infty(M)\) compactly supported 
the following formula holds
\begin{align}
\label{eq:identity-biharmonic-noncompact}
&\int_M\eta^2\big(\lambda(m-4)-\operatorname{Scal}^M\big)|\tau(\phi)|^2dv_g+
4\int_M\eta^2\operatorname{Ric}^M(e_i,e_j)\langle\bar\nabla d\phi(e_i,e_j),\tau(\phi)\rangle dv_g \\
\nonumber=&
-4\int_M\langle d\phi(\operatorname{Ric}^M(\nabla\eta^2)),\tau(\phi)\rangle dv_g 
+4\lambda\int_M\langle d\phi(\nabla\eta^2),\tau(\phi)\rangle dv_g 
-\int_M\langle\nabla\eta^2,\nabla f\rangle|\tau(\phi)|^2 dv_g \\
&\nonumber+2\int_M(\Delta\eta^2)\langle d\phi(\nabla f),\tau(\phi)\rangle dv_g 
+4\int_M(\nabla_{e_i}\eta^2)\nabla_{e_j}f\langle\bar\nabla d\phi(e_i,e_j),\tau(\phi)\rangle dv_g,
\end{align}
where \(\{e_i\},i=1,\ldots,m\) represents an orthonormal basis of \(TM\).
\end{Lem}

\begin{proof}
By assumption \(\phi\) is a smooth biharmonic map such that the stress-energy tensor
associated with biharmonic maps \eqref{stress-energy-bienery} is divergence free.
As \(\eta\) has compact support we use integration by parts to deduce
\begin{equation*}
0=-\int_M\eta^2\langle\nabla f,\operatorname{div} S_2\rangle dv_g
=\int_M\eta^2\langle\nabla^2f,S_2\rangle dv_g
+\int_M S_2(\nabla f,\nabla\eta^2) dv_g.
\end{equation*}

In addition, we calculate
\begin{align*}
\int_M\eta^2\langle g,S_2\rangle dv_g=&\int_M\eta^2\tr_gS_2 dv_g \\
=&\int_M\eta^2\big(m(\frac{1}{2}|\tau(\phi)|^2+\langle d\phi,\bar\nabla\tau(\phi)\rangle)
-2\langle d\phi,\bar\nabla\tau(\phi)\rangle\big) dv_g\\
=&\big(2-\frac{m}{2}\big)\int_M\eta^2|\tau(\phi)|^2 dv_g
+(2-m)\int_M\langle d\phi(\nabla\eta^2),\tau(\phi)\rangle dv_g.
\end{align*}
Moreover, a similar calculation shows
\begin{align*}
\int_M\eta^2\langle\operatorname{Ric}^M,S_2\rangle dv_g
=&-\frac{1}{2}\int_M\eta^2\operatorname{Scal}^M|\tau(\phi)|^2 dv_g
+2\int_M\eta^2\operatorname{Ric}^M(e_i,e_j)\langle\bar\nabla d\phi(e_i,e_j),\tau(\phi)\rangle dv_g\\
&+2\int_M\eta^2\langle d\phi\big(\underbrace{\operatorname{div}\operatorname{Ric}^M-\frac{1}{2}\nabla \operatorname{Scal}^M}_{=0}\big),\tau(\phi)\rangle dv_g \\
&-\int_M \operatorname{Scal}^M
\langle d\phi(\nabla\eta^2),\tau(\phi)\rangle dv_g
+2\int_M\langle d\phi(\operatorname{Ric}^M(\nabla\eta^2)),\tau(\phi)\rangle dv_g.
\end{align*}

As a next step we manipulate the term \(S_2(\nabla f,\nabla\eta^2)\).
To this end, let \(\{e_i\},i=1,\ldots,m\) be an orthonormal basis of \(TM\) that satisfies 
\(\nabla_{e_i}e_j=0,i,j=1,\ldots,m\) at a fixed point \(p\in M\).
Using the definition of the stress-energy tensor \eqref{stress-energy-bienery}
we find
\begin{align*}
S_2(\nabla f,\nabla\eta^2)
=&\langle\nabla\eta^2,\nabla f\rangle\big(\frac{1}{2}|\tau(\phi)|^2
+\langle d\phi,\bar\nabla\tau(\phi)\rangle\big)
-(\nabla_{e_i}\eta^2)\nabla_{e_j}f\langle d\phi(e_i),\bar\nabla_{e_j}\tau(\phi)\rangle \\
&-(\nabla_{e_j}\eta^2)\nabla_{e_i}f\langle d\phi(e_i),\bar\nabla_{e_j}\tau(\phi)\rangle.
\end{align*}

A direct calculation using integration by parts yields
\begin{align*}
\int_M\langle\nabla\eta^2,\nabla f\rangle\big(\frac{1}{2}|\tau(\phi)|^2
+\langle d\phi,\bar\nabla\tau(\phi)\rangle\big) dv_g
=&-\frac{1}{2}\int_M\langle\nabla\eta^2,\nabla f\rangle|\tau(\phi)|^2 dv_g \\
&-\int_M(\nabla_{e_j}\nabla_{e_i}\eta^2)\nabla_{e_i}f\langle d\phi(e_j),\tau(\phi)\rangle dv_g\\
&-\int_M(\nabla_{e_i}\eta^2)\nabla_{e_j}\nabla_{e_i}f\langle d\phi(e_j),\tau(\phi)\rangle dv_g.
\end{align*}

Moreover, using integration by parts once more we find
\begin{align*}
\int_M(\nabla_{e_i}\eta^2)\nabla_{e_j} f\langle d\phi(e_i),\bar\nabla_{e_j}\tau(\phi)\rangle dv_g 
=&-\int_M(\nabla_{e_j}\nabla_{e_i}\eta^2)\nabla_{e_j} f\langle d\phi(e_i),\tau(\phi)\rangle dv_g \\
&-\int_M(\nabla_{e_i}\eta^2)\Delta f\langle d\phi(e_i),\tau(\phi)\rangle dv_g \\
&-\int_M(\nabla_{e_i}\eta^2)\nabla_{e_j} f\langle\bar\nabla d\phi(e_i,e_j),\tau(\phi)\rangle dv_g
\end{align*}
and 
\begin{align*}
\int_M(\nabla_{e_j}\eta^2)\nabla_{e_i} f\langle d\phi(e_i),\bar\nabla_{e_j}\tau(\phi)\rangle dv_g 
=&-\int_M(\Delta\eta^2)\nabla_{e_i} f\langle d\phi(e_i),\tau(\phi)\rangle dv_g \\
&-\int_M(\nabla_{e_j}\eta^2)\nabla_{e_j}\nabla_{e_i}f
\langle d\phi(e_i),\tau(\phi)\rangle dv_g \\
&-\int_M(\nabla_{e_j}\eta^2)\nabla_{e_i} f\langle\bar\nabla d\phi(e_i,e_j),\tau(\phi)\rangle dv_g.
\end{align*}

Adding up both contributions yields
\begin{align*}
&\int_M(\nabla_{e_i}\eta^2)\nabla_{e_j}f\langle d\phi(e_i),\bar\nabla_{e_j}\tau(\phi)\rangle dv_g 
+\int_M(\nabla_{e_j}\eta^2)\nabla_{e_i}f\langle d\phi(e_i),\bar\nabla_{e_j}\tau(\phi)\rangle dv_g \\
=&-\int_M(\nabla_{e_j}\nabla_{e_i}\eta^2)\nabla_{e_j}f\langle d\phi(e_i),\tau(\phi)\rangle dv_g
-\int_M(\nabla_{e_i}\eta^2)\Delta f\langle d\phi(e_i),\tau(\phi)\rangle dv_g \\
&-2\int_M(\nabla_{e_i}\eta^2)\nabla_{e_j}f\langle\bar\nabla d\phi(e_i,e_j),\tau(\phi)\rangle dv_g
-\int_M(\Delta\eta^2)\nabla_{e_i}f\langle d\phi(e_i),\tau(\phi)\rangle dv_g \\
&-\int_M(\nabla_{e_j}\eta^2)\nabla_{e_j}\nabla_{e_i}f\langle d\phi(e_i),\tau(\phi)\rangle dv_g.
\end{align*}

Combining the previous equations we get
\begin{align*}
\int_M S_2(\nabla f,\nabla\eta^2)~dv_g=&
-\frac{1}{2}\int_M\langle\nabla\eta^2,\nabla f\rangle|\tau(\phi)|^2 dv_g 
+\int_M(\Delta\eta^2)\langle d\phi(\nabla f),\tau(\phi)\rangle dv_g \\
\nonumber&+\int_M\Delta f\langle d\phi(\nabla\eta^2),\tau(\phi)\rangle dv_g  
+2\int_M(\nabla_{e_i}\eta^2)\nabla_{e_j}f\langle\bar\nabla d\phi(e_i,e_j),\tau(\phi)\rangle dv_g\\
=&-\frac{1}{2}\int_M\langle\nabla\eta^2,\nabla f\rangle|\tau(\phi)|^2 dv_g 
+\int_M(\Delta\eta^2)\langle d\phi(\nabla f),\tau(\phi)\rangle dv_g \\
\nonumber&+\int_M(m\lambda-\operatorname{Scal}^M)\langle d\phi(\nabla\eta^2),\tau(\phi)\rangle dv_g \\
&+2\int_M(\nabla_{e_i}\eta^2)\nabla_{e_j}f\langle\bar\nabla d\phi(e_i,e_j),\tau(\phi)\rangle dv_g,
\end{align*}
where we used the equation for a gradient Ricci soliton \eqref{eq:soliton} 
in the second step.

The result follows from adding up the different contributions.
\end{proof}

\begin{proof}[Proof of Theorem \ref{thm:liouville-biharmonic}]
Again, let  \(0\leq\eta\leq 1\) on \(M\) be such that
\begin{equation*}
\eta(x)=1\textrm{ for } x\in B_R(x_0),\qquad \eta(x)=0\textrm{ for } x\in B_{2R}(x_0),\qquad |\nabla^q\eta|\leq\frac{C}{R^q}\textrm{ for } x\in M,
\end{equation*}
where \(B_R(x_0)\) denotes the geodesic ball around the point \(x_0\) with radius \(R\) and \(q=1,2\).
Moreover, let \(\{e_i\},i=1,\ldots,m\) be an orthonormal basis of \(TM\) that satisfies 
\(\nabla_{e_i}e_j=0,i,j=1,\ldots,m\) at a fixed point \(p\in M\).

In order to prove the result we estimate the terms on 
the right hand side of \eqref{eq:identity-biharmonic-noncompact}.
Making use of the finite energy assumption \eqref{eq:biharmonic-finite-energy}
it is easy to infer
\begin{align*}
\int_M\langle d\phi(\operatorname{Ric}^M(\nabla\eta^2)),\tau(\phi)\rangle dv_g
\leq&\frac{C}{R}\int_M(|\bar\nabla d\phi|^2+|d\phi|^2)dv_g \to 0, \\
\int_M\langle d\phi(\nabla\eta^2),\tau(\phi)\rangle dv_g 
\leq&\frac{C}{R}\int_M(|\bar\nabla d\phi|^2+|d\phi|^2)dv_g 
\to 0,  \\
\int_M\langle\nabla\eta^2,\nabla f\rangle|\tau(\phi)|^2 dv_g 
\leq&\frac{C}{R}\int_M|\bar\nabla d\phi|^2dv_g 
\to 0, \\
\nonumber\int_M(\Delta\eta^2)\langle d\phi(\nabla f),\tau(\phi)\rangle dv_g 
\leq&\frac{C}{R^2}\int_M(|\bar\nabla d\phi|^2+|d\phi|^2)dv_g 
\to 0, \\
\nonumber\int_M(\nabla_{e_i}\eta^2)\nabla_{e_j}f\langle\bar\nabla d\phi(e_i,e_j),\tau(\phi)\rangle 
\leq&\frac{C}{R}\int_M|\bar\nabla d\phi|^2 dv_g\to 0
\end{align*}
as \(R\to\infty\).
Note that we used the inequality \(|\tau(\phi)|\leq\sqrt{m}|\bar\nabla d\phi|\).
This completes the proof.
\end{proof}

\begin{proof}[Proof of Theorem \ref{cor:cigar-biharmonic}]
Recall that the cigar soliton is a steady gradient Ricci soliton
meaning that it is a solution of \eqref{eq:soliton} with \(\lambda=0\).
Since we are on a two-dimensional domain we also have
\(2\operatorname{Ric}^M=\operatorname{Scal}^Mg\).

In the case of a steady gradient Ricci soliton we always
get a pointwise bound on \(\nabla f\), see Lemma \ref{lem:bound-nabla-f}
and it is also straightforward to see that we have an upper bound 
on the Ricci curvature of the cigar soliton.
Hence, we may apply Theorem \ref{thm:liouville-harmonic} and 
inserting the above data into \eqref{identity-biharmonic} yields
\begin{equation*}
\int_{\R^2}\operatorname{Scal}^M|\tau(\phi)|^2dv_g\leq 0.
\end{equation*}
As the cigar soliton has positive scalar curvature we can deduce \(\tau(\phi)=0\)
yielding the claim.
\end{proof}

\section{Some remarks about harmonic and biharmonic maps from Yamabe solitons}
Another class of Riemannian manifolds that is closely connected to gradient Ricci solitons
are manifolds that admit a concircular vector field.
In this case we have
\begin{equation}
\label{eq:concircular}
\nabla^2F=\varphi g
\end{equation}
where both \(F,\varphi\) are smooth functions on \(M\).
Such kinds of manifolds have been intensively studied by Tashiro \cite{MR174022}.

A special case of \eqref{eq:concircular} is the equation for
a gradient Yamabe soliton given by
\begin{equation}
\label{eq:yamabe}
\nabla^2F=(\operatorname{Scal}^M-\rho)g
\end{equation}
where \(\operatorname{Scal}^M\) represents the scalar curvature of \(M\)
and \(\rho\in\R\). 
These solitons arise as self-similar solutions of the Yamabe flow and have 
been classified in \cite{MR3008413} and \cite{MR2989649}.

As in the case of gradient Ricci solitons, a Yamabe soliton is called 
\emph{steady} if \(\rho=0\),
\emph{shrinking} if \(\rho>0\) and \emph{expanding} if \(\rho<0\).

We have seen that in the case of a steady gradient Ricci soliton we automatically
get a pointwise bound on \(\nabla f\), see Lemma \ref{lem:bound-nabla-f}.
There does not seem to be a corresponding result for gradient Yamabe solitons,
see \cite{MR3573991} for more details.

However, it is straightforward to prove the following 
\begin{Satz}
\label{thm:liouville-yamabe}
Assume that \((M,g,F)\) is a complete non-compact gradient Yamabe soliton with 
 \(|\nabla F|<\infty\).
\begin{enumerate}
 \item Let \(\phi\colon M\to N\) be a smooth harmonic map with finite energy, that is
\begin{equation*}
\int_M|d\phi|^2dv_g<\infty.
\end{equation*}
Then the following inequality holds
\begin{equation}
\label{inequality-harmonic-yamabe}
(m-2)\int_M\big(\operatorname{Scal}^M-\rho\big)|d\phi|^2dv_g
\leq 0.
\end{equation}
\item Let \(\phi\colon M\to N\) be a smooth biharmonic map with finite energy, that is
\begin{equation*}
\int_M(|d\phi|^2+|\bar\nabla d\phi|^2)dv_g<\infty.
\end{equation*}
Then the following inequality holds
\begin{equation}
\label{inequality-biharmonic-yamabe}
(m-4)\int_M\big(\operatorname{Scal}^M-\rho\big)|\tau(\phi)|^2dv_g
\leq 0.
\end{equation}
\end{enumerate}
\end{Satz}

One can of course deduce vanishing results from Theorem \ref{thm:liouville-yamabe} but we do not further investigate this matter here.

\bibliographystyle{wileyNJD-AMS}

\begin{thebibliography}{10}

\bibitem{MR655417}
P.~Baird and J.~Eells.
\newblock A conservation law for harmonic maps.
\newblock In {\em Geometry {S}ymposium, {U}trecht 1980 ({U}trecht, 1980)},
  volume 894 of {\em Lecture Notes in Math.}, pages 1--25. Springer, Berlin-New
  York, 1981.

\bibitem{MR2604617}
Paul Baird, Ali Fardoun, and Seddik Ouakkas.
\newblock Liouville-type theorems for biharmonic maps between {R}iemannian
  manifolds.
\newblock {\em Adv. Calc. Var.}, 3(1):49--68, 2010.

\bibitem{MR3340295}
Jacob Bernstein and Thomas Mettler.
\newblock Two-dimensional gradient {R}icci solitons revisited.
\newblock {\em Int. Math. Res. Not. IMRN}, (1):78--98, 2015.

\bibitem{MR4106647}
V.~Branding, S.~Montaldo, C.~Oniciuc, and A.~Ratto.
\newblock Higher order energy functionals.
\newblock {\em Adv. Math.}, 370:107236, 60, 2020.

\bibitem{MR3834926}
Volker Branding.
\newblock A {L}iouville-type theorem for biharmonic maps between complete
  {R}iemannian manifolds with small energies.
\newblock {\em Arch. Math. (Basel)}, 111(3):329--336, 2018.

\bibitem{MR4142862}
Volker Branding.
\newblock Some analytic results on interpolating sesqui-harmonic maps.
\newblock {\em Ann. Mat. Pura Appl. (4)}, 199(5):2039--2059, 2020.

\bibitem{MR4007262}
Volker Branding.
\newblock The stress-energy tensor for polyharmonic maps.
\newblock {\em Nonlinear Anal.}, 190:111616, 17, 2020.

\bibitem{MR4040175}
Volker Branding and Yong Luo.
\newblock A nonexistence theorem for proper biharmonic maps into general
  {R}iemannian manifolds.
\newblock {\em J. Geom. Phys.}, 148:103557, 9, 2020.

\bibitem{MR3008413}
Huai-Dong Cao, Xiaofeng Sun, and Yingying Zhang.
\newblock On the structure of gradient {Y}amabe solitons.
\newblock {\em Math. Res. Lett.}, 19(4):767--774, 2012.

\bibitem{MR2989649}
Giovanni Catino, Carlo Mantegazza, and Lorenzo Mazzieri.
\newblock On the global structure of conformal gradient solitons with
  nonnegative {R}icci tensor.
\newblock {\em Commun. Contemp. Math.}, 14(6):1250045, 12, 2012.

\bibitem{MR2520796}
Bing-Long Chen.
\newblock Strong uniqueness of the {R}icci flow.
\newblock {\em J. Differential Geom.}, 82(2):363--382, 2009.

\bibitem{MR2448435}
Manolo Eminenti, Gabriele La~Nave, and Carlo Mantegazza.
\newblock Ricci solitons: the equation point of view.
\newblock {\em Manuscripta Math.}, 127(3):345--367, 2008.

\bibitem{MR3406677}
Hong~Xin Guo, Robert Philipowski, and Anton Thalmaier.
\newblock On gradient solitons of the {R}icci-harmonic flow.
\newblock {\em Acta Math. Sin. (Engl. Ser.)}, 31(11):1798--1804, 2015.

\bibitem{MR891928}
Guo~Ying Jiang.
\newblock The conservation law for {$2$}-harmonic maps between {R}iemannian
  manifolds.
\newblock {\em Acta Math. Sinica}, 30(2):220--225, 1987.

\bibitem{MR2395125}
E.~Loubeau, S.~Montaldo, and C.~Oniciuc.
\newblock The stress-energy tensor for biharmonic maps.
\newblock {\em Math. Z.}, 259(3):503--524, 2008.

\bibitem{MR2301373}
S.~Montaldo and C.~Oniciuc.
\newblock A short survey on biharmonic maps between {R}iemannian manifolds.
\newblock {\em Rev. Un. Mat. Argentina}, 47(2):1--22 (2007), 2006.

\bibitem{MR2961788}
Reto M\"{u}ller.
\newblock Ricci flow coupled with harmonic map flow.
\newblock {\em Ann. Sci. \'{E}c. Norm. Sup\'{e}r. (4)}, 45(1):101--142, 2012.

\bibitem{MR3023848}
Ovidiu Munteanu and Natasa Sesum.
\newblock On gradient {R}icci solitons.
\newblock {\em J. Geom. Anal.}, 23(2):539--561, 2013.

\bibitem{MR4265170}
Ye-Lin Ou and Bang-Yen Chen.
\newblock {\em Biharmonic submanifolds and biharmonic maps in {R}iemannian
  geometry}.
\newblock World Scientific Publishing Co. Pte. Ltd., Hackensack, NJ, [2020]
  \copyright 2020.

\bibitem{MR2507581}
Peter Petersen and William Wylie.
\newblock Rigidity of gradient {R}icci solitons.
\newblock {\em Pacific J. Math.}, 241(2):329--345, 2009.

\bibitem{MR2818729}
Stefano Pigola, Michele Rimoldi, and Alberto~G. Setti.
\newblock Remarks on non-compact gradient {R}icci solitons.
\newblock {\em Math. Z.}, 268(3-4):777--790, 2011.

\bibitem{ram}
Daniel Ramos.
\newblock Gradient {R}icci solitons on surfaces.
\newblock {\em arXiv preprint arXiv:1304.6391}, 2013.

\bibitem{MR438388}
Richard Schoen and Shing~Tung Yau.
\newblock Harmonic maps and the topology of stable hypersurfaces and manifolds
  with non-negative {R}icci curvature.
\newblock {\em Comment. Math. Helv.}, 51(3):333--341, 1976.

\bibitem{MR654088}
H.~C.~J. Sealey.
\newblock Some conditions ensuring the vanishing of harmonic differential forms
  with applications to harmonic maps and {Y}ang-{M}ills theory.
\newblock {\em Math. Proc. Cambridge Philos. Soc.}, 91(3):441--452, 1982.

\bibitem{MR174022}
Yoshihiro Tashiro.
\newblock Complete {R}iemannian manifolds and some vector fields.
\newblock {\em Trans. Amer. Math. Soc.}, 117:251--275, 1965.

\bibitem{MR2265040}
Peter Topping.
\newblock {\em Lectures on the {R}icci flow}, volume 325 of {\em London
  Mathematical Society Lecture Note Series}.
\newblock Cambridge University Press, Cambridge, 2006.

\bibitem{MR1391729}
Yuanlong Xin.
\newblock {\em Geometry of harmonic maps}, volume~23 of {\em Progress in
  Nonlinear Differential Equations and their Applications}.
\newblock Birkh\"{a}user Boston, Inc., Boston, MA, 1996.

\bibitem{MR3573991}
Fei Yang and Liangdi Zhang.
\newblock Geometry of gradient {Y}amabe solitons.
\newblock {\em Ann. Global Anal. Geom.}, 50(4):367--379, 2016.

\bibitem{MR3668864}
Gabjin Yun.
\newblock {$L^2$} harmonic forms on gradient shrinking {R}icci solitons.
\newblock {\em J. Korean Math. Soc.}, 54(4):1189--1208, 2017.

\end{thebibliography}

\end{document}